\documentclass{amsart}

\usepackage[utf8x]{inputenc}

\usepackage{amssymb}
\usepackage{amsmath}
\usepackage{amscd}
\usepackage{units}

\usepackage[english]{babel}

\usepackage{tikz}
\usetikzlibrary{arrows}

\newtheorem{theorem}{Theorem}[section]

\newtheorem{proposition}[theorem]{Proposition}

\newtheorem{definition}[theorem]{Definition}

\newtheorem{conjecture}[theorem]{Conjecture}

\usepackage{color}

\title{Three-dimensional Catalan numbers and product-coproduct prographs}

\author{Nicolas Borie}

\address{Université Paris-Est, LIGM (UMR 8049), UPEM, CNRS, ENPC,
  ESIEE, F-77454, Marne-la-Vallée, France.}


\begin{document}

\maketitle

\begin{abstract}
  We present the new combinatorial class of product-coproduct
  prographs which are planar assemblies of two types of operators:
  products having two inputs and a single output and coproducts having
  a single input and two outputs. We show that such graphs are enumerated
  by the $3$-dimensional Catalan numbers. We present some
  combinatorial bijections positioning product-coproduct prographs as key
  objects to probe families of objects enumerated by the $3$-dimensional
  Catalan numbers.
\end{abstract}

\section{Introduction}

This story begins with computer explorations. Recall that planar
rooted binary trees are planar structures freely generated by a formal
operator having one output (a single father per node) and two inputs
(left and right children); this single operator can be viewed as a
non-associative product. Now, if we add a non-coassociative coproduct
with a single input and two outputs, what will be the set of planar
structures built from these two formal operators? To bound this
problem, we restrict the enumeration to structures having globally a
single input and a single output, thus each enumerated element will
contain as many products as coproducts.

Such an element is called a
\emph{prograph} with one input and one output.

\begin{figure}[h]
  \centering
\begin{tikzpicture}[scale=0.55]
  \tikzstyle{prod}=[fill,draw,rectangle,minimum size=5pt,inner sep=1pt]
  \tikzstyle{cop}=[fill,draw,circle,minimum size=6pt,inner sep=1pt]

  \draw (0, 0.5) -- (0, 0);
  \draw (0,0) -- (-0.5,-0.5) -- (0, -1);
  \draw (0,0) -- (0.5,-0.5) -- (0, -1);
  \draw (0, -1) -- (0, -2);
  \draw (0,-2) -- (-0.5,-2.5) -- (0, -3);
  \draw (0,-2) -- (0.5,-2.5) -- (0, -3);
  \draw (0, -3) -- (0, -3.5);

  \draw (0,0) node[prod] (p2) {$~$};
  \draw (0,-1) node[cop] (c2) {$~$};
  \draw (0,-2) node[prod] (p1) {$~$};
  \draw (0,-3) node[cop] (c1) {$~$};

  \draw (4,0.5) -- (4,0);
  \draw (3,-1) -- (4,0) -- (5, -1);
  \draw (2,-2) -- (3,-1) -- (4, -2);
  \draw (5,-1) -- (4,-2);
  \draw (2,-2) -- (3,-3) -- (4,-2);
  \draw (3,-3) -- (3,-3.5);

  \draw (4,0) node[prod] (p2) {$~$};
  \draw (3,-1) node[prod] (p1) {$~$};
  \draw (4,-2) node[cop] (c2) {$~$};
  \draw (3,-3) node[cop] (c1) {$~$};

  \draw (8,0.5) -- (8,0);
  \draw (8,0) -- (9.5,-1.5) -- (8, -3);
  \draw (8,0) -- (6.5,-1.5) -- (8, -3);
  \draw (7,-1) -- (7.5,-1.5) -- (7, -2);
  \draw (8,-3) -- (8,-3.5);

  \draw (8,0) node[prod] (p2) {$~$};
  \draw (7,-1) node[prod] (p1) {$~$};
  \draw (7,-2) node[cop] (c2) {$~$};
  \draw (8,-3) node[cop] (c1) {$~$};

  \draw (12,0.5) -- (12,0);
  \draw (12,0) -- (13.5,-1.5) -- (12, -3);
  \draw (12,0) -- (10.5,-1.5) -- (12, -3);
  \draw (13,-1) -- (12.5,-1.5) -- (13, -2);
  \draw (12,-3) -- (12,-3.5);

  \draw (12,0) node[prod] (p2) {$~$};
  \draw (13,-1) node[prod] (p1) {$~$};
  \draw (13,-2) node[cop] (c2) {$~$};
  \draw (12,-3) node[cop] (c1) {$~$};

  \draw (16,0.5) -- (16,0);
  \draw (16,0) -- (15,-1) -- (17, -3);
  \draw (16,0) -- (18,-2) -- (17, -3);
  \draw (17,-1) -- (16,-2);
  \draw (17,-3) -- (17,-3.5);

  \draw (16,0) node[prod] (p2) {$~$};
  \draw (17,-1) node[prod] (p1) {$~$};
  \draw (16,-2) node[cop] (c2) {$~$};
  \draw (17,-3) node[cop] (c1) {$~$};

\end{tikzpicture}
\caption{The five prographs with two coproducts, two products,
  and having a single input and a single output.}~\label{Fig1}
\end{figure}
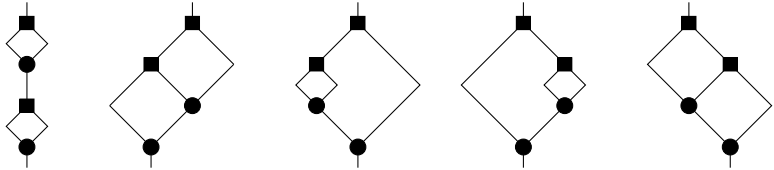

We implemented a Sage~\cite{sage} program enumerating
prographs with $n$ products and $n$ coproducts and the first
values were $1$, $1$, $5$, $42$, $462$, $6006$, $87516$. This is
the beginning of Sequence A005789 of the
OEIS~\cite{Sloane}, the $3$-dimensional Catalan numbers.
Figure~\ref{Fig1} displays the 
$5$ prographs containing two coproducts (circles) and two products
(squares). We will see how one can show that prographs containing
$n$ products and coproducts are counted by the $n^{th}$
$3$-dimensional Catalan numbers. Moreover, we will show that prographs
are relevant structures on which a lot of combinatorics can be expressed.

In this paper, we present four combinatorial classes
enumerated by the $3$-dimensional Catalan numbers and we show how their
combinatorics are related to product-coproduct prographs. In the
following section, we recall the definitions and some properties of
$3$-dimensional Catalan numbers, standard Young tableaux with three rows and
product-coproduct prographs.
In Section~\ref{section_edges}, we present a well-chosen labeling of
edges of prographs which can be extended to an isomorphism of Hopf
algebras. Section~\ref{section_up_down} presents up-down permutations
of even size that avoid the pattern $(1234)$. Trying to establish another
combinatorial bijection with prographs, we also present in the same
section a bijection between permutations avoiding $(123)$ and binary
trees (to the best of our knowledge, this bijection is not present in
the literature). Finally, in Section~\ref{section_boxes}, we show that
a labeling of operators of prographs gives another Hopf isomorphism
with weighted Dyck paths having some constraints; such paths are
also Laguerre histories. This will allow us to build a
partial combinatorial bijection between prographs and up-down
permutations of even size that avoid the pattern $(1234)$.

\section{Preliminaries}


\begin{definition}
  For $n$ a non-negative integer, the $n^{th}$ $3$-dimensional Catalan
  number counts the number of paths from $(0,0,0)$ to $(n,n,n)$ using
  steps $(+1,0,0)$, $(0,+1,0)$ and $(0,0,+1)$ such that each point
  $(x,y,z)$ on the path satisfies $x \geqslant y \geqslant z$.
\end{definition}

It is obvious that the $n^{th}$ $3$-dimensional Catalan
number also counts the number of standard tableaux of shape
$(n,n,n)$: from a $3$-dimensional path, read from left to right, 
assign a number from $1$ to $3n$ to each step and fill the tableau
by inserting labels of step $(+1,0,0)$ on the first row, $(0,+1,0)$ on the
second row and $(0,0,+1)$ on the third row. The following example
uses the French convention for tableaux.

\begin{small}
  \begin{equation}
    \begin{array}{ccccccc}
      (0,0,0) & \rightarrow (1,0,0) & \rightarrow (1,1,0) & \rightarrow
      (2,1,0) & \rightarrow (2,1,1) & \rightarrow (2,2,1) & \rightarrow
      (2,2,2) \\
              & +~(1,0,0) & +~(0,1,0) & +~(1,0,0) & +~(0,0,1) & +~(0,1,0) & +~(0,0,1) \\
              & 1 & 2 & 3 & 4 & 5 & 6 \\
    \end{array}
    \begin{array}{|c|c|} \hline
      4 & 6 \\ \hline
      2 & 5 \\ \hline
      1 & 3 \\ \hline
    \end{array}
\end{equation}
\end{small}

We denote by $ST_{\langle n^3 \rangle}$ the set of standard Young
tableaux of shape $(n,n,n)$.  The hook-length formula for standard
tableaux gives the following nice expression.

\begin{proposition}
  The $n^{th}$ $3$-dimensional Catalan number $C^{(3)}_n$ is given by:
  \begin{equation}
    |ST_{\langle n^3 \rangle}| = \frac{2\cdot(3n)!}{n! \cdot (n+1)! \cdot (n+2)!}.
  \end{equation}
\end{proposition}


Although it is complicated to describe in general, the Schützenberger
involution of rectangular standard Young tableaux can be obtained with an
easy algorithm.

\begin{proposition}[See~\cite{stan_schutz, Pon_Wang} for details]
  The Schützenberger involution $S$ has a simple description on rectangular
  standard Young tableaux: it consists in reversing the alphabet
  $\{1, 2, \dots , n \}$ and rotating the tableaux by $180^{\circ}$.
\end{proposition}

Formalized by category theory~\cite{Mac_Lane}, PROs are often viewed as the
natural generalization of operads. These are still planar assemblies
of formal operators, but now, formal operators have not necessarily a
single output. The product over objects remains grafting, however
objects are not trees anymore but graphs.

Free PROs over a finite set of
generators constitute a nice example of PRO to apprehend this
object. Given a finite set of operators $G$, the free PRO generated by
$G$ is the set of all finite planar graphs (called \emph{prographs}) freely
built using elements in $G$ (each operator can appear several
times). As ``free'' means also formal, a generator inside $G$ can
just be described by its number of inputs and outputs.

\begin{definition}
A \emph{product-coproduct prograph} is a connected directed graph having
a single input and a single output and composed with two types of
nodes: coproduct nodes having a single input and two outputs
and product nodes having two inputs and a single output.
\end{definition}

As each coproduct introduces a new output and each product
suppresses an output, a product-coproduct prograph must
contain as many products as coproducts. For $n$ a non-negative integer,
we will denote by $PC(n)$ the set of product coproduct prographs
containing $n$ coproducts and $n$ products.

Let $V$ be a formal module, let
$\Delta : V \rightarrow V \otimes V$ be a coproduct and
$\cdot : V \otimes V \rightarrow V$ a product. Encoding the operator
\begin{tikzpicture}[scale=0.6]
  \tikzset{
    boxi/.style={fill,draw,rectangle,minimum size=5pt,inner sep=1pt}
  }
  \tikzstyle{prod}=[fill,draw,rectangle,minimum size=5pt,inner sep=1pt]
  \tikzstyle{cop}=[fill,draw,circle,minimum size=6pt,inner sep=1pt]

  \draw (0,0) -- (0,-0.4);
  \draw (0,0) -- (0.3,0.3);
  \draw (0,0) -- (-0.3,0.3);

  \draw (0,0) node[cop] (c1) {$~$};
\end{tikzpicture} having a single entry and two outputs with $\Delta$ and encoding
\begin{tikzpicture}[scale=0.6]
  \tikzset{
    boxi/.style={fill,draw,rectangle,minimum size=5pt,inner sep=1pt}
  }
  \tikzstyle{prod}=[fill,draw,rectangle,minimum size=5pt,inner sep=1pt]
  \tikzstyle{cop}=[fill,draw,circle,minimum size=6pt,inner sep=1pt]

  \draw (0,0) -- (0,0.4);
  \draw (0,0) -- (0.3,-0.3);
  \draw (0,0) -- (-0.3,-0.3);

  \draw (0,0) node[prod] (c1) {$~$};
\end{tikzpicture} with the product $\cdot$, we can associate
an algebraic expression with each prograph. The associated expression
models a map from $V$ to $V$ and is a composition of layers
which are mainly tensor products of some $\Delta$, $\cdot$ and
the identity map $Id$. Without any operator, the only product-coproduct prograph is the
unique empty prograph. Doing nothing (or do not change anything)
corresponds to the map $Id$.

There is a single prograph with one coproduct and one product: $\cdot \circ \Delta$.

With two coproducts and two products, we get the five following
expressions associated with the five prographs of Figure~\ref{Fig1} from
left to right:
\begin{equation}
  \begin{array}{c}
    \cdot \circ \Delta \circ \cdot \circ \Delta, \qquad
    \cdot \circ (\cdot \otimes Id) \circ (Id \otimes \Delta) \circ \Delta, \qquad 
    \cdot \circ (\cdot \otimes Id) \circ (\Delta \otimes Id) \circ \Delta, \\ 
    \cdot \circ (Id \otimes \cdot) \circ (Id \otimes \Delta) \circ \Delta, \qquad
    \cdot \circ (Id \otimes \cdot) \circ (\Delta \otimes Id) \circ \Delta.
    \end{array}
\end{equation}


\begin{definition}
  Rotating by $180^{\circ}$ naturally defines an involution $S$ on
  prographs which we will call the Schützenberger involution.
\end{definition}
We will see in Proposition~\ref{prop_schutz} that it is equivalent to
the classical Schützenberger involution, hence the name.

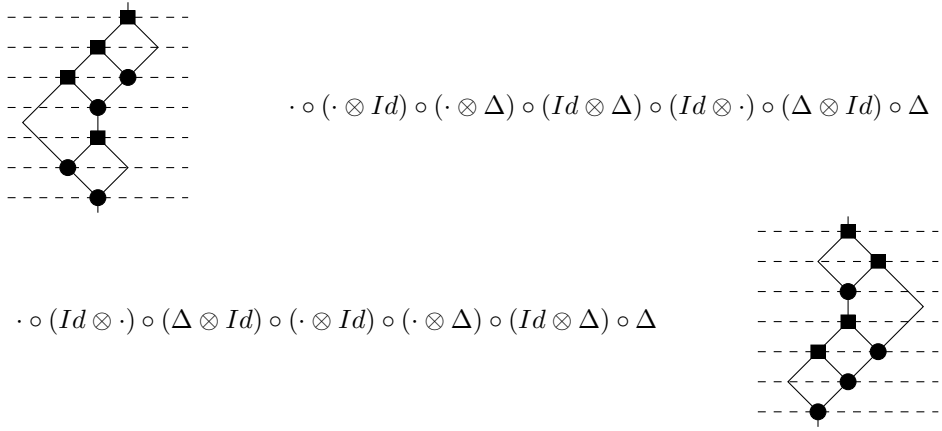
\begin{figure}[h]
  \centering
\begin{tikzpicture}[scale=0.4]
  \tikzset{
    boxi/.style={fill,draw,rectangle,minimum size=5pt,inner sep=1pt}
  }
  \tikzstyle{prod}=[fill,draw,rectangle,minimum size=5pt,inner sep=1pt]
  \tikzstyle{cop}=[fill,draw,circle,minimum size=6pt,inner sep=1pt]

  \draw (0,0) -- (0,-0.5);
  \draw (0,0) -- (1,1) -- (0,2);
  \draw (0,0) -- (-2.5,2.5) -- (1,6);
  \draw (0,2) -- (0,3);
  \draw (-1,1) -- (0,2);
  \draw (0,3) -- (-1,4);
  \draw (1,4) -- (0,5);
  \draw (0,3) -- (2,5) -- (1,6);
  \draw (1,6) -- (1,6.5);

  \draw (0,0) node[cop] (c1) {$~$};
  \draw (-1,1) node[cop] (c2) {$~$};
  \draw (0,2) node[prod] (p1) {$~$};
  \draw (0,3) node[cop] (c3) {$~$};  
  \draw (-1,4) node[prod] (p2) {$~$};
  \draw (1,4) node[cop] (c4) {$~$};
  \draw (0,5) node[prod] (p3) {$~$};
  \draw (1,6) node[prod] (p4) {$~$};

  \draw[dashed] (-3,0) -- (3,0);
  \draw[dashed] (-3,1) -- (3,1);
  \draw[dashed] (-3,2) -- (3,2);
  \draw[dashed] (-3,3) -- (3,3);
  \draw[dashed] (-3,4) -- (3,4);
  \draw[dashed] (-3,5) -- (3,5);
  \draw[dashed] (-3,6) -- (3,6);
  
  \draw (17,3) node (exp) {$\cdot \circ (\cdot \otimes Id) \circ (\cdot \otimes \Delta) \circ (Id \otimes \Delta) \circ (Id \otimes \cdot) \circ (\Delta \otimes Id) \circ \Delta $};

\end{tikzpicture} \\
\begin{tikzpicture}[scale=0.4]
  \tikzset{
    boxi/.style={fill,draw,rectangle,minimum size=5pt,inner sep=1pt}
  }
  \tikzstyle{prod}=[fill,draw,rectangle,minimum size=5pt,inner sep=1pt]
  \tikzstyle{cop}=[fill,draw,circle,minimum size=6pt,inner sep=1pt]

  \draw (5,0) -- (5,-0.5);
  \draw (5,0) -- (4,1) -- (6,3);
  \draw (5,0) -- (8.5,3.5) -- (6,6);
  \draw (5,2) -- (6,1);
  \draw (6,3) -- (6,4);
  \draw (6,3) -- (7,2);
  \draw (6,4) -- (7,5);
  \draw (6,4) -- (5,5) -- (6,6);
  \draw (6,6) -- (6,6.5);

  \draw (5,0) node[cop] (c1) {$~$};
  \draw (6,1) node[cop] (c2) {$~$};
  \draw (7,2) node[cop] (c3) {$~$};
  \draw (5,2) node[prod] (p1) {$~$};
  \draw (6,3) node[prod] (p2) {$~$};
  \draw (6,4) node[cop] (c4) {$~$};  
  \draw (7,5) node[prod] (p3) {$~$};
  \draw (6,6) node[prod] (p4) {$~$};

  \draw[dashed] (3,0) -- (9,0);
  \draw[dashed] (3,1) -- (9,1);
  \draw[dashed] (3,2) -- (9,2);
  \draw[dashed] (3,3) -- (9,3);
  \draw[dashed] (3,4) -- (9,4);
  \draw[dashed] (3,5) -- (9,5);
  \draw[dashed] (3,6) -- (9,6);

  \draw (-11,3) node (exp) {$\cdot \circ (Id \otimes \cdot) \circ (\Delta \otimes Id) \circ (\cdot \otimes Id) \circ (\cdot \otimes \Delta) \circ (Id \otimes \Delta) \circ \Delta $};
  
\end{tikzpicture}
  
\caption{A prograph, its image by the Schützenberger involution and their algebraic expressions.}~\label{algexpr_prograph}

\end{figure}

The example of Figure~\ref{algexpr_prograph} is a little bigger one with
seven layers on the graph, and thus, the expression describing the
product coproduct prograph has seven blocks of operators.

Reading the expression from left to right (respectively from right to
left) corresponds to scanning the prograph from top to bottom
(respectively from bottom to top). On the algebraic expression, the
Schützenberger involution consists in switching coproducts $\Delta$
and products $\cdot$, and reversing the obtained expression.

As our product-coproduct prographs have a single input and a single
output, we can assemble prographs by stacking them: the output of the
first one grafted with the input of the second one. This operation
defines a product making the disjoint union
$PC := \bigcup_{n \in \mathbb{N}} PC(n)$ a monoid.
The algebra of this monoid, coupled with the proper coproduct, turns
out to be a Hopf algebra.

\section{Labeling edges of product-coproduct prographs}~\label{section_edges}

A well-chosen labeling of the edges of prographs gives a first
bijection, which strongly motivated our investigations on these
objects.

\begin{figure}[h]
  \centering
  \scalebox{0.7}{
\begin{tikzpicture}
  \tikzset{
    boxi/.style={fill,draw,rectangle,minimum size=5pt,inner sep=1pt}
  }
  \tikzstyle{prod}=[fill,draw,rectangle,minimum size=5pt,inner sep=1pt]
  \tikzstyle{cop}=[fill,draw,circle,minimum size=6pt,inner sep=1pt]

  \draw (0,0) -- (0,-0.5);
  \draw (0,0) -- (1,1) -- (0,2);
  \draw (0,0) -- (-2.5,2.5) -- (1,6);
  \draw (0,2) -- (0,3);
  \draw (-1,1) -- (0,2);
  \draw (0,3) -- (-1,4);
  \draw (1,4) -- (0,5);
  \draw (0,3) -- (2,5) -- (1,6);
  \draw (1,6) -- (1,6.5);

  \draw (0,0) node[cop] (c1) {$~$};
  \draw (-1,1) node[cop] (c2) {$~$};
  \draw (0,2) node[prod] (p1) {$~$};
  \draw (0,3) node[cop] (c3) {$~$};  
  \draw (-1,4) node[prod] (p2) {$~$};
  \draw (1,4) node[cop] (c4) {$~$};
  \draw (0,5) node[prod] (p3) {$~$};
  \draw (1,6) node[prod] (p4) {$~$};

  \draw (0,-0.25) node[left] (1) {$1$};
  \draw (-0.5,0.4) node[left] (2) {$2$};
  \draw (-2.5,2.5) node[left] (3) {$3$};
  \draw (-0.5,1.65) node[left] (4) {$4$};
  \draw (0.9,1) node[left] (5) {$5$};
  \draw (0,2.5) node[left] (6) {$6$};
  \draw (-0.5,3.35) node[left] (7) {$7$};
  \draw (-0.5,4.65) node[left] (8) {$8$};
  \draw (0.5,3.65) node[left] (9) {$9$};
  \draw (0.5,4.4) node[left] (10) {$10$};
  \draw (0.5,5.65) node[left] (11) {$11$};
  \draw (1.85,5) node[left] (12) {$12$};

  \draw (3,3) node (expr) {$Sch.$};
  \draw (3,2.5) node (expr) {$\longleftrightarrow$};  

  \draw (5,0) -- (5,-0.5);
  \draw (5,0) -- (4,1) -- (6,3);
  \draw (5,0) -- (8.5,3.5) -- (6,6);
  \draw (5,2) -- (6,1);
  \draw (6,3) -- (6,4);
  \draw (6,3) -- (7,2);
  \draw (6,4) -- (7,5);
  \draw (6,4) -- (5,5) -- (6,6);
  \draw (6,6) -- (6,6.5);

  \draw (5,0) node[cop] (c1) {$~$};
  \draw (6,1) node[cop] (c2) {$~$};
  \draw (7,2) node[cop] (c3) {$~$};
  \draw (5,2) node[prod] (p1) {$~$};
  \draw (6,3) node[prod] (p2) {$~$};
  \draw (6,4) node[cop] (c4) {$~$};  
  \draw (7,5) node[prod] (p3) {$~$};
  \draw (6,6) node[prod] (p4) {$~$};

  \draw (5,-0.25) node[left] (1) {$1$};
  \draw (4,1) node[left] (2) {$2$};
  \draw (5.5,0.65) node[left] (3) {$3$};
  \draw (5.5,1.35) node[left] (4) {$4$};
  \draw (5.5,2.65) node[left] (5) {$5$};
  \draw (6.5,1.65) node[left] (6) {$6$};
  \draw (6.5,2.35) node[left] (7) {$7$};
  \draw (6,3.5) node[left] (8) {$8$};
  \draw (5,5) node[left] (9) {$9$};
  \draw (6.5,4.65) node[left] (10) {$10$};
  \draw (8.35,3.5) node[left] (11) {$11$};
  \draw (6.5,5.35) node[left] (12) {$12$};
\end{tikzpicture}}
\caption{Labeling of the edges of a prograph and its rotated version.}~\label{labeledges}
\end{figure}
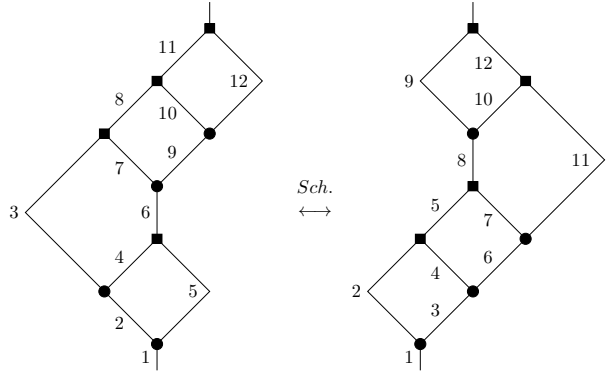

\begin{theorem}
  The set $PC(n)$ of product-coproduct prographs of size $n$ has as
  cardinality the $n^{th}$ 3-dimensional Catalan number.
  \begin{proof}
    We build a bijection $le$ (labeling edges)
    \begin{displaymath}
      le : \displaystyle\bigcup_{n \in \mathbb{N}} PC(n) \rightarrow
      \displaystyle\bigcup_{n \in \mathbb{N}} ST_{\langle n^3 \rangle}
    \end{displaymath}
    and its inverse using a depth-left first search numbering of
    wires on prographs (See Figure~\ref{labeledges}). The first row
    of the tableau will contain the labels of the inputs of coproducts,
    the second row will contain the left inputs of products, and the
    third row will contains right inputs. 
  \end{proof}
\end{theorem}

\begin{proposition}~\label{prop_schutz}
  The bijection $le$ preserves the Schützenberger involution. That is, we have:
  \begin{displaymath}
    S(p) = le^{-1}(  S(  le(p))).
  \end{displaymath}
  \begin{proof}
    The depth-left first labeling forms a covering path from bottom to
    top of the prograph (with indices from $1$ to $3n$). On the
    rotated version, we take the same path in its reverse way,
    labeling with $i$ instead $3n-i$.
  \end{proof}
\end{proposition}

The filling of Young tableaux of prographs of Figure~\ref{labeledges}
and the Schützenberger involution give

\begin{equation}
  \begin{array}{|c|c|c|c|}\hline
    5 & 7 & 10 & 12 \\ \hline
    3 & 4 & 8 & 11 \\ \hline
    1 & 2 & 6 & 9 \\ \hline
  \end{array}
  \quad
  \begin{array}{c}
    \text{reverse} \\
    \text{alphabet} \\
    \rightarrow
  \end{array}
  \quad
  \begin{array}{|c|c|c|c|}\hline
    8 & 6 & 3 & 1 \\ \hline
    10 & 9 & 5 & 2 \\ \hline
    12 & 11 & 7 & 4 \\ \hline
  \end{array}
  \quad
  \begin{array}{c}
    \text{rotation} \\
    \text{by } 180^{\circ} \\
    \rightarrow
  \end{array}
  \quad
  \begin{array}{|c|c|c|c|}\hline
    4 & 7 & 11 & 12 \\ \hline
    2 & 5 & 9 & 10 \\ \hline
    1 & 3 & 6 & 8 \\ \hline
  \end{array}~.
\end{equation}

The shifted concatenation product $\bullet$ gives a monoid structure
over three-row standard Young tableaux. For instance,

\begin{equation}
  \begin{array}{|c|c|c|} \hline
    5 & 8 & 9 \\ \hline
    3 & 4 & 7 \\ \hline
    1 & 2 & 6 \\ \hline
  \end{array}
  \bullet
  \begin{array}{|c|c|} \hline
    5 & 6 \\ \hline
    2 & 4 \\ \hline
    1 & 3 \\ \hline
  \end{array}
  =
  \begin{array}{|c|c|c|c|c|} \hline
    5 & 8 & 9 & 14 & 15 \\ \hline
    3 & 4 & 7 & 11 & 13 \\ \hline
    1 & 2 & 6 & 10 & 12 \\ \hline
  \end{array}~.
\end{equation}

\begin{proposition}
   For $\mathbb{K}$ a field, the bijection $le$ extended by linearity
   on the monoid algebra $\mathbb{K}[\bigcup_{n \in \mathbb{N}} PC(n)]$
   with values inside the monoid algebra
   $\mathbb{K}[\bigcup_{n \in \mathbb{N}} ST_{\langle n^3 \rangle}]$
   becomes an isomorphism of graded Hopf algebras.
\end{proposition}

\section{Up-down permutations of $2n$ avoiding $(1234)$}~\label{section_up_down}

An up-down permutation of $2n$ avoiding $(1234)$ is a
permutation of size $2n$ whose descents set is $\{2, 4, 6, \dots \}$
with no four letters that form an increasing subsequence. We denote by $A_{2n}(1234)$
the set formed by all these permutations.

Here are the $42$ up-down permutations of size $6$ avoiding $(1234)$:
\begin{equation}
  A_{6}(1234) = \left\{
  \begin{array}{c}
    563412, 562413, 562314, 561423, 561324, 463512, 462513, \\
    462315, 461523, 461325, 453612, 452613, 452316, 451623, \\
    451326, 364512, 362514, 362415, 361524, 361425, 354612, \\
    352614, 352416, 351624, 351426, 342615, 341625, 264513, \\
    263514, 261534, 261435, 254613, 253614, 251634, 251436, \\
    243615, 241635, 164523, 163524, 154623, 153624, 143625 
  \end{array}
  \right\}.
\end{equation}

Lewis proved that the $n^{th}$ $3$-dimensional Catalan numbers
count the cardinality of $A_{2n}(1234)$ and gives in~\cite{Lewis1,
  Lewis2} two bijections between up-down permutations of $2n$
avoiding $(1234)$ and standard Young tableaux of shape $(n,n,n)$. However,
using these two bijections, we did not manage to prove that up-down
permutations of $2n$ avoiding $(1234)$ deploy the same combinatorics
as standard Young tableaux or prographs. At first glance, these results
appear to us mainly as counting results.

On permutations, we also have the classical Schützenberger involution
(and we will once more denote it by $S$) which consists in reversing
the alphabet, then reversing the reading direction. For example,
$S(631278594) = 615238974$. As the Schützenberger involution
preserves appearance and avoidance of patterns, $S$ stabilizes the set
of up-down permutations of $2n$ avoiding $(1234)$.

\begin{definition}
  We define a shifted concatenation product $\bullet$ on
  $\bigcup_{n \in \mathbb{N}} A_{2n}(1234)$ as
  \begin{equation}
    \begin{array}{ccccl}
      \displaystyle\bigcup_{n \in \mathbb{N}} A_{2n}(1234) & \otimes & \displaystyle\bigcup_{n \in \mathbb{N}} A_{2n}(1234) & \longrightarrow & \displaystyle\bigcup_{n \in \mathbb{N}} A_{2n}(1234) \\
      (\qquad \sigma & , & \tau \qquad) & \longmapsto & (\operatorname{shift}_{length(\sigma)}(\tau)) \cdot \sigma \\
    \end{array}
  \end{equation}
\end{definition}

Here are some examples:
\begin{equation}
  \begin{array}{c}
  12 \bullet 12 = 3412, \\
  2143 \bullet 1324 = 57682143, \\
  12^{\bullet n} = (2n-1)(2n)(2n-3)(2n -2) \cdots 563412. \\
  \end{array}
\end{equation}

Not only having the same cardinality, we think that up-down
permutations avoiding $(1234)$ present the same combinatorics as
prographs. This is formulated in the following conjecture.

\begin{conjecture}~\label{conj1234}
  There exists a bijection between up-down permutations of $2n$
  avoiding $(1234)$ and prographs of $PC(n)$ preserving the
  Schützenberger involution that can be extended to an isomorphism
  of Hopf algebras.
\end{conjecture}

For $\sigma$ an up-down permutation of $2n$ avoiding $(1234)$, we
will denote by $Peaks(\sigma)$ the subsequence of values at even
positions and $Vals(\sigma)$ the subsequence of values at odd positions
inside $\sigma$.

\begin{proposition}~\label{prop_cond}
  For $n$ a non-negative integer and $\sigma$ a permutation of size
  $2n$, $\sigma$ is an up-down permutation of $2n$ avoiding $(1234)$
  if and only if the following four conditions are verified:
  \begin{itemize}
    \item the sequence $Peaks(\sigma)$ avoids $(123)$, and
    \item the sequence $Vals(\sigma)$ avoids $(123)$, and
    \item each value of $Peaks(\sigma)$ lower than a valley $k$ appears
      to the right of $k$ in $\sigma$, and
    \item if a valley $k$ has a lower valley to its left, all peak
      values greater than $k$ to its right must be ordered in
      $\sigma$ decreasingly.
  \end{itemize}
  \begin{proof}
    By exhaustion of all possible positions of values that
    would form a 1234 pattern.
  \end{proof}
\end{proposition}

Proposition~\ref{prop_cond} gathers conditions not very
handy for describing up-down permutations of $2n$ avoiding
$(1234)$. However, since permutations
avoiding $(123)$ are counted by (classical) Catalan numbers and thus,
are in bijection with binary trees, this proposition presents the
remaining conditions we will need to build special product-coproduct
prographs from a pair of binary trees (the second being reversed
over the first one).

Let us now build a bijection between binary trees and permutations
avoiding $(123)$ compatible with the depth-left labeling algorithm.
Let $\sigma$ be a permutation of size $n-1$ avoiding $(123)$. The
possible positions of a new value $n$ to be inserted in $\sigma$ such
that it still avoids $(123)$ are constrained. Let $\tau$ the maximal
prefix of $\sigma$ whose values are decreasing. If $\sigma$ begins by
a rise, $\tau$ contains only the first value of $\sigma$. If
$\sigma$ is entirely decreasing, then $\tau = \sigma$. The possible
positions to insert $n$ in $\sigma$ are before $\tau$, just after
$\tau$ or inside $\tau$. By inserting $n$ farther, we would
get a new permutation $\tau \sigma_1 n \sigma_2$ where
$\tau \sigma_1 \sigma_2 = \sigma$ and $\sigma_1 \neq \epsilon$. Such a permutation
would contain for sure a pattern $(123)$ where the smallest value can
be in $\tau$, the middle one in $\sigma_1$ and the value $n$ for the
greatest one. At each insertion of the largest value, the number of
values after the first rise (or zero if the permutation is entirely
decreasing) is a non-decreasing statistic bounded by the size of the
permutation minus one. This gives a way to identify a permutation
avoiding $(123)$ with a non-decreasing parking function (we mean here
a non-decreasing function from $\{1, \dots , n\}$ to $\{0, \dots ,n-1\}$
such that $f(i) \leqslant i-1$).

On the other side, when one labels a planar tree (drawn from bottom to
top) from the root with a depth-left first algorithm, at each insertion (or
new label), the number of insertion positions left free on the left
is non-decreasing and bounded by the number of nodes, and therefore,
forms a non-decreasing parking function. The construction of the tree
associated with the permutation $958732641$ is presented Figure~\ref{bij123_tree}.

\begin{figure}
  \centering
  \scalebox{0.79}{
  \begin{tabular}{|cc|cc|cc|} \hline
    \begin{tabular}{c}
      1 \\
      \small{(0)} \\
    \end{tabular} &
    \begin{tikzpicture}[scale=0.65]
      \tikzstyle{cop}=[draw,circle,minimum size=6pt,inner sep=1pt]
      
      \draw (0,0) node[cop] (1) {$1$};
    \end{tikzpicture} &
    \begin{tabular}{c}
      21 \\
      \small{(0,0)} \\
    \end{tabular} &
    \begin{tikzpicture}[scale=0.65]
      \tikzstyle{cop}=[draw,circle,minimum size=6pt,inner sep=1pt]
      
      \draw (0,0) node[cop] (1) {$1$};
      \draw (-1,1) node[cop] (2) {$2$};
      \draw (1) edge (2);
    \end{tikzpicture} &
    \begin{tabular}{c}
      321 \\
      \small{(0,0,0)} \\
    \end{tabular} &
    \begin{tikzpicture}[scale=0.65]
      \tikzstyle{cop}=[draw,circle,minimum size=6pt,inner sep=1pt]
      
      \draw (0,0) node[cop] (1) {$1$};
      \draw (-1,1) node[cop] (2) {$2$};
      \draw (-2,2) node[cop] (3) {$3$};
      \draw (0,2.35) node (x) {$~$};
      \draw (1) edge (2);
      \draw (2) edge (3);
    \end{tikzpicture} \\ \hline
    \begin{tabular}{c}
      32\underline{\textbf{41}} \\
      \small{(0,0,0,2)} \\
    \end{tabular} &
    \begin{tikzpicture}[scale=0.65]
      \tikzstyle{cop}=[draw,circle,minimum size=6pt,inner sep=1pt]
      
      \draw (0,0) node[cop] (1) {$1$};
      \draw (-1,1) node[cop] (2) {$2$};
      \draw (-2,2) node[cop] (3) {$3$};
      \draw (0,2) node[cop] (4) {$4$};
      \draw (-2.35, 2.35) node (x1) {$\bullet$};
      \draw (-1.65, 2.35) node (x2) {$\bullet$};
      \draw (1) edge (2);
      \draw (2) edge (3);
      \draw (2) edge (4);
    \end{tikzpicture} &
    \begin{tabular}{c}
      532\underline{\textbf{41}} \\
      \small{(0,0,0,2,2)} \\
    \end{tabular} &
    \begin{tikzpicture}[scale=0.65]
      \tikzstyle{cop}=[draw,circle,minimum size=6pt,inner sep=1pt]
      
      \draw (0,0) node[cop] (1) {$1$};
      \draw (-1,1) node[cop] (2) {$2$};
      \draw (-2,2) node[cop] (3) {$3$};
      \draw (0,2) node[cop] (4) {$4$};
      \draw (-2.35, 2.35) node (x1) {$\bullet$};
      \draw (-1.65, 2.35) node (x2) {$\bullet$};
      \draw (-1,3) node[cop] (5) {$5$};
      \draw (1) edge (2);
      \draw (2) edge (3);
      \draw (2) edge (4);
      \draw (4) edge (5);
    \end{tikzpicture} &
    \begin{tabular}{c}
      532\underline{\textbf{641}} \\
      \small{(0,0,0,2,2,3)} \\
    \end{tabular} &
    \begin{tikzpicture}[scale=0.65]
      \tikzstyle{cop}=[draw,circle,minimum size=6pt,inner sep=1pt]
      
      \draw (0,0) node[cop] (1) {$1$};
      \draw (-1,1) node[cop] (2) {$2$};
      \draw (-2,2) node[cop] (3) {$3$};
      \draw (0,2) node[cop] (4) {$4$};
      \draw (-2.35, 2.35) node (x1) {$\bullet$};
      \draw (-1.65, 2.35) node (x2) {$\bullet$};
      \draw (-1,3) node[cop] (5) {$5$};
      \draw (0,4) node[cop] (6) {$6$};
      \draw (-1.35, 3.35) node (x3) {$\bullet$};
      \draw (0.35, 4.35) node (x) {$~$};
      \draw (1) edge (2);
      \draw (2) edge (3);
      \draw (2) edge (4);
      \draw (4) edge (5);
      \draw (5) edge (6);
    \end{tikzpicture} \\ \hline
    \begin{tabular}{c}
      5\underline{\textbf{732641}} \\
      \small{(0,0,0,2,2,3,6)} \\
    \end{tabular} &
    \begin{tikzpicture}[scale=0.65]
      \tikzstyle{cop}=[draw,circle,minimum size=6pt,inner sep=1pt]
      
      \draw (0,0) node[cop] (1) {$1$};
      \draw (-1,1) node[cop] (2) {$2$};
      \draw (-2,2) node[cop] (3) {$3$};
      \draw (0,2) node[cop] (4) {$4$};
      \draw (-2.35, 2.35) node (x1) {$\bullet$};
      \draw (-1.65, 2.35) node (x2) {$\bullet$};
      \draw (-1,3) node[cop] (5) {$5$};
      \draw (0,4) node[cop] (6) {$6$};
      \draw (-1.35, 3.35) node (x3) {$\bullet$};
      \draw (1,1) node[cop] (7) {$7$};
      \draw (-0.35, 4.35) node (x4) {$\bullet$};
      \draw (0.35, 4.35) node (x5) {$\bullet$};
      \draw (0.35, 2.35) node (x6) {$\bullet$};
      \draw (1) edge (2);
      \draw (2) edge (3);
      \draw (2) edge (4);
      \draw (4) edge (5);
      \draw (5) edge (6);
      \draw (1) edge (7);
    \end{tikzpicture} &
    \begin{tabular}{c}
      5\underline{\textbf{8732641}} \\
      \small{(0,0,0,2,2,3,6,7)} \\
    \end{tabular} &
    \begin{tikzpicture}[scale=0.65]
      \tikzstyle{cop}=[draw,circle,minimum size=6pt,inner sep=1pt]
      
      \draw (0,0) node[cop] (1) {$1$};
      \draw (-1,1) node[cop] (2) {$2$};
      \draw (-2,2) node[cop] (3) {$3$};
      \draw (0,2) node[cop] (4) {$4$};
      \draw (-2.35, 2.35) node (x1) {$\bullet$};
      \draw (-1.65, 2.35) node (x2) {$\bullet$};
      \draw (-1,3) node[cop] (5) {$5$};
      \draw (0,4) node[cop] (6) {$6$};
      \draw (-1.35, 3.35) node (x3) {$\bullet$};
      \draw (1,1) node[cop] (7) {$7$};
      \draw (-0.35, 4.35) node (x4) {$\bullet$};
      \draw (0.35, 4.35) node (x5) {$\bullet$};
      \draw (0.35, 2.35) node (x6) {$\bullet$};
      \draw (2,2) node[cop] (8) {$8$};
      \draw (0.65, 1.35) node (x7) {$\bullet$};
      \draw (1) edge (2);
      \draw (2) edge (3);
      \draw (2) edge (4);
      \draw (4) edge (5);
      \draw (5) edge (6);
      \draw (1) edge (7);
      \draw (7) edge (8);
    \end{tikzpicture} &
        \begin{tabular}{c}
      95\underline{\textbf{8732641}} \\
      \small{(0,0,0,2,2,3,6,7,7)} \\
    \end{tabular} &
    \begin{tikzpicture}[scale=0.65]
      \tikzstyle{cop}=[draw,circle,minimum size=6pt,inner sep=1pt]
      
      \draw (0,0) node[cop] (1) {$1$};
      \draw (-1,1) node[cop] (2) {$2$};
      \draw (-2,2) node[cop] (3) {$3$};
      \draw (0,2) node[cop] (4) {$4$};
      \draw (-2.35, 2.35) node (x1) {$\bullet$};
      \draw (-1.65, 2.35) node (x2) {$\bullet$};
      \draw (-1,3) node[cop] (5) {$5$};
      \draw (0,4) node[cop] (6) {$6$};
      \draw (-1.35, 3.35) node (x3) {$\bullet$};
      \draw (1,1) node[cop] (7) {$7$};
      \draw (-0.35, 4.35) node (x4) {$\bullet$};
      \draw (0.35, 4.35) node (x5) {$\bullet$};
      \draw (0.35, 2.35) node (x6) {$\bullet$};
      \draw (2,2) node[cop] (8) {$8$};
      \draw (0.65, 1.35) node (x7) {$\bullet$};
      \draw (1,3) node[cop] (9) {$9$};
      \draw (1) edge (2);
      \draw (2) edge (3);
      \draw (2) edge (4);
      \draw (4) edge (5);
      \draw (5) edge (6);
      \draw (1) edge (7);
      \draw (7) edge (8);
      \draw (8) edge (9);
    \end{tikzpicture} \\ \hline
  \end{tabular}}
  \caption{Insertion algorithms for permutations avoiding $(123)$,
    non-decreasing parking functions and binary trees labeled by
    depth-left first traversal.}~\label{bij123_tree}
\end{figure}

Figure~\ref{trees_14} presents our bijection between the $14$
permutations of size $4$ avoiding $(123)$ and the $14$ binary trees
having $4$ nodes.

\begin{figure}
  \centering
  \begin{tabular}{ccc|ccc|ccc}
    4321 & (0,0,0,0) &
    \begin{tikzpicture}[scale=0.3]
      \draw (0,0) node (1) {$\bullet$};
      \draw (-1,1) node (1) {$\bullet$};
      \draw (-2,2) node (1) {$\bullet$};
      \draw (-3,3) node (1) {$\bullet$};
      \draw (0,0) -- (-3.5, 3.5);
      \draw (0,0) -- (0,-0.5);
      \draw (0,0) -- (0.5,0.5);
      \draw (-1,1) -- (-0.5,1.5);
      \draw (-2,2) -- (-1.5,2.5);
      \draw (-3,3) -- (-2.5,3.5);
    \end{tikzpicture}
    &
    3214 & (0,0,0,1) &
    \begin{tikzpicture}[scale=0.3]
      \draw (0,0) node (1) {$\bullet$};
      \draw (-1,1) node (1) {$\bullet$};
      \draw (-2,2) node (1) {$\bullet$};
      \draw (-1,3) node (1) {$\bullet$};
      \draw (-1,3) -- (-1.5, 3.5);
      \draw (0,0) -- (0,-0.5);
      \draw (0,0) -- (0.5,0.5);
      \draw (0,0) -- (-2.5,2.5);
      \draw (-1,1) -- (-0.5,1.5);
      \draw (-2,2) -- (-0.5,3.5);
    \end{tikzpicture}
    &
    4213 & (0,0,1,1) &
    \begin{tikzpicture}[scale=0.3]
      \draw (0,0) node (1) {$\bullet$};
      \draw (-1,1) node (1) {$\bullet$};
      \draw (0,2) node (1) {$\bullet$};
      \draw (-1,3) node (1) {$\bullet$};
      \draw (0,2) -- (-1.5, 3.5);
      \draw (0,0) -- (0,-0.5);
      \draw (0,0) -- (0.5,0.5);
      \draw (0,0) -- (-1.5,1.5);
      \draw (-1,1) -- (0.5,2.5);
      \draw (-1,3) -- (-0.5,3.5);
    \end{tikzpicture}
    \\ \hline
    4312 & (0,1,1,1) &
    \begin{tikzpicture}[scale=0.3]
      \draw (0,0) node (1) {$\bullet$};
      \draw (1,1) node (1) {$\bullet$};
      \draw (0,2) node (1) {$\bullet$};
      \draw (-1,3) node (1) {$\bullet$};
      \draw (1,1) -- (-1.5, 3.5);
      \draw (0,0) -- (0,-0.5);
      \draw (0,0) -- (1.5,1.5);
      \draw (0,0) -- (-0.5,0.5);
      \draw (0,2) -- (0.5,2.5);
      \draw (-1,3) -- (-0.5,3.5);
    \end{tikzpicture}
    &
    3241 & (0,0,0,2) &
    \begin{tikzpicture}[scale=0.3]
      \draw (0,0) node (1) {$\bullet$};
      \draw (-1,1) node (1) {$\bullet$};
      \draw (-2,2) node (1) {$\bullet$};
      \draw (0,2) node (1) {$\bullet$};
      \draw (0,2) -- (-0.5, 2.5);
      \draw (0,0) -- (0,-0.5);
      \draw (0,0) -- (0.5,0.5);
      \draw (0,0) -- (-2.5,2.5);
      \draw (-1,1) -- (0.5,2.5);
      \draw (-2,2) -- (-1.5,2.5);
    \end{tikzpicture}
    &
    2143 & (0,0,1,2) &
    \begin{tikzpicture}[scale=0.3]
      \draw (0,0) node (1) {$\bullet$};
      \draw (-1,1) node (1) {$\bullet$};
      \draw (0,2) node (1) {$\bullet$};
      \draw (1,3) node (1) {$\bullet$};
      \draw (0,2) -- (-0.5, 2.5);
      \draw (0,0) -- (0,-0.5);
      \draw (0,0) -- (0.5,0.5);
      \draw (0,0) -- (-1.5,1.5);
      \draw (-1,1) -- (1.5,3.5);
      \draw (1,3) -- (0.5,3.5);
    \end{tikzpicture}
    \\ \hline
    3142 & (0,1,1,2) &
    \begin{tikzpicture}[scale=0.3]
      \draw (0,0) node (1) {$\bullet$};
      \draw (1,1) node (1) {$\bullet$};
      \draw (0,2) node (1) {$\bullet$};
      \draw (1,3) node (1) {$\bullet$};
      \draw (1,1) -- (-0.5, 2.5);
      \draw (0,0) -- (0,-0.5);
      \draw (0,0) -- (1.5,1.5);
      \draw (0,0) -- (-0.5,0.5);
      \draw (0,2) -- (1.5,3.5);
      \draw (1,3) -- (0.5,3.5);
    \end{tikzpicture}
    &
    4231 & (0,0,2,2) &
    \begin{tikzpicture}[scale=0.3]
      \draw (0,0) node (1) {$\bullet$};
      \draw (-1,1) node (1) {$\bullet$};
      \draw (0,2) node (1) {$\bullet$};
      \draw (1,1) node (1) {$\bullet$};
      \draw (0,0) -- (0,-0.5);
      \draw (0,0) -- (1.5,1.5);
      \draw (0,0) -- (-1.5,1.5);
      \draw (1,1) -- (-0.5,2.5);
      \draw (0,2) -- (0.5,2.5);
      \draw (-1,1) -- (-0.5,1.5);
    \end{tikzpicture}
    &
    4132 & (0,1,2,2) &
        \begin{tikzpicture}[scale=0.3]
      \draw (0,0) node (1) {$\bullet$};
      \draw (1,3) node (1) {$\bullet$};
      \draw (2,2) node (1) {$\bullet$};
      \draw (1,1) node (1) {$\bullet$};
      \draw (0,0) -- (0,-0.5);
      \draw (0,0) -- (2.5,2.5);
      \draw (0,0) -- (-0.5,0.5);
      \draw (1,1) -- (0.5, 1.5);
      \draw (2,2) -- (0.5, 3.5);
      \draw (1,3) -- (1.5,3.5);
    \end{tikzpicture}
    \\ \hline
    3421 & (0,0,0,3) &
    \begin{tikzpicture}[scale=0.3]
      \draw (0,0) node (1) {$\bullet$};
      \draw (-1,1) node (1) {$\bullet$};
      \draw (-2,2) node (1) {$\bullet$};
      \draw (1,1) node (1) {$\bullet$};
      \draw (0,0) -- (-2.5, 2.5);
      \draw (0,0) -- (0,-0.5);
      \draw (0,0) -- (1.5,1.5);
      \draw (-1,1) -- (-0.5,1.5);
      \draw (-2,2) -- (-1.5,2.5);
      \draw (1,1) -- (0.5,1.5);
    \end{tikzpicture}    
    &
    2413 & (0,0,1,3) &
    \begin{tikzpicture}[scale=0.3]
      \draw (0,0) node (1) {$\bullet$};
      \draw (-1,1) node (1) {$\bullet$};
      \draw (0,2) node (1) {$\bullet$};
      \draw (1,1) node (1) {$\bullet$};
      \draw (0,0) -- (0,-0.5);
      \draw (0,0) -- (1.5,1.5);
      \draw (0,0) -- (-1.5,1.5);
      \draw (1,1) -- (0.5,1.5);
      \draw (-1,1) -- (0.5,2.5);
      \draw (0,2) -- (-0.5,2.5);
    \end{tikzpicture}    
    &
    3412 & (0,1,1,3) &
    \begin{tikzpicture}[scale=0.3]
      \draw (0,0) node (1) {$\bullet$};
      \draw (2,2) node (1) {$\bullet$};
      \draw (0,2) node (1) {$\bullet$};
      \draw (1,1) node (1) {$\bullet$};
      \draw (0,0) -- (0,-0.5);
      \draw (0,0) -- (2.5,2.5);
      \draw (0,0) -- (-0.5,0.5);
      \draw (1,1) -- (-0.5,2.5);
      \draw (0,2) -- (0.5,2.5);
      \draw (2,2) -- (1.5,2.5);
    \end{tikzpicture}
    \\ \hline
    2431 & (0,0,2,3) &
    \begin{tikzpicture}[scale=0.3]
      \draw (0,0) node (1) {$\bullet$};
      \draw (-1,1) node (1) {$\bullet$};
      \draw (2,2) node (1) {$\bullet$};
      \draw (1,1) node (1) {$\bullet$};
      \draw (0,0) -- (-1.5, 1.5);
      \draw (0,0) -- (0,-0.5);
      \draw (0,0) -- (2.5,2.5);
      \draw (-1,1) -- (-0.5,1.5);
      \draw (2,2) -- (1.5,2.5);
      \draw (1,1) -- (0.5,1.5);
    \end{tikzpicture}
    &
    1432 & (0,1,2,3) &
    \begin{tikzpicture}[scale=0.3]
      \draw (0,0) node (1) {$\bullet$};
      \draw (3,3) node (1) {$\bullet$};
      \draw (2,2) node (1) {$\bullet$};
      \draw (1,1) node (1) {$\bullet$};
      \draw (0,0) -- (-0.5, 0.5);
      \draw (0,0) -- (0,-0.5);
      \draw (0,0) -- (3.5,3.5);
      \draw (3,3) -- (2.5,3.5);
      \draw (2,2) -- (1.5,2.5);
      \draw (1,1) -- (0.5,1.5);
    \end{tikzpicture}
    \\
  \end{tabular}
  \caption{The $14$ permutations of size $4$ avoiding $(123)$, their
    corresponding non-decreasing parking functions and binary trees.}~\label{trees_14}
\end{figure}

\section{Labeling boxes of product-coproduct prographs}~\label{section_boxes}

After having labeled the edges of prographs and recovered the three-row
standard Young tableaux, it seems natural to investigate what we
obtain when we label operators (boxes in the prographs). For $n$ a
non-negative integer, a prograph of $PC(n)$ contains $n$ coproducts
and $n$ products, therefore, the labels will run from $1$ up to
$2n$. We will still use depth-left first algorithm to label operators
of prographs since it preserves the Schützenberger involution. 

\begin{figure}[h]
  \centering
\begin{tikzpicture}[scale=1]
  \tikzset{
    boxi/.style={fill,draw,rectangle,minimum size=5pt,inner sep=1pt}
  }
  \tikzstyle{prod}=[draw,rectangle,minimum size=7pt,inner sep=2pt]
  \tikzstyle{cop}=[draw,circle,minimum size=6pt,inner sep=1pt]
  \tikzstyle{virt}=[minimum size=0pt,inner sep=0pt]  

  \draw (0,0) node[cop] (c1) {$1$};
  \draw (-1,1) node[cop] (c2) {$2$};
  \draw (0,2) node[prod] (p1) {$3$};
  \draw (0,3) node[cop] (c3) {$4$};  
  \draw (-1,4) node[prod] (p2) {$5$};
  \draw (1,4) node[cop] (c4) {$6$};
  \draw (0,5) node[prod] (p3) {$7$};
  \draw (1,6) node[prod] (p4) {$8$};

  \draw (1,1) node[virt] (e1) {};
  \draw (-2.5,2.5) node[virt] (e2) {};
  \draw (2,5) node[virt] (e3) {};
  \draw (0,-0.5) node[virt] (e4) {};
  \draw (1,6.5) node[virt] (e5) {};
  
  \draw (c1) edge (c2);
  \draw (c1) edge (e1);
  \draw (e1) edge (p1);
  \draw (c2) edge (e2);
  \draw (e2) edge (p2);
  \draw (c2) edge (p1);
  \draw (p1) edge (c3);
  \draw (c3) edge (p2);
  \draw (c3) edge (c4);
  \draw (p3) edge (p4);
  \draw (p2) edge (p3);
  \draw (c4) edge (p3);
  \draw (c4) edge (e3);
  \draw (e3) edge (p4);
  \draw (p4) edge (e5);
  \draw (c1) edge (e4);

  \draw (5,0) node[cop] (c1) {$1$};
  \draw (6,1) node[cop] (c2) {$2$};
  \draw (7,2) node[cop] (c3) {$4$};
  \draw (5,2) node[prod] (p1) {$3$};
  \draw (6,3) node[prod] (p2) {$5$};
  \draw (6,4) node[cop] (c4) {$6$};  
  \draw (7,5) node[prod] (p3) {$7$};
  \draw (6,6) node[prod] (p4) {$8$};

  \draw (4,1) node[virt] (e1) {};
  \draw (8.5,3.5) node[virt] (e2) {};
  \draw (5,5) node[virt] (e3) {};
  \draw (5,-0.5) node[virt] (e4) {};
  \draw (6,6.5) node[virt] (e5) {};

  \draw (c1) edge (c2);
  \draw (c1) edge (e1);
  \draw (e1) edge (p1);
  \draw (c3) edge (e2);
  \draw (e2) edge (p3);
  \draw (c2) edge (p1);
  \draw (p1) edge (p2);
  \draw (c3) edge (p2);
  \draw (p2) edge (c4);
  \draw (p3) edge (p4);
  \draw (c2) edge (c3);
  \draw (c4) edge (p3);
  \draw (c4) edge (e3);
  \draw (e3) edge (p4);
  \draw (p4) edge (e5);
  \draw (c1) edge (e4);
\end{tikzpicture}
\caption{Labeling the operators of a prograph and its rotation.}~\label{labelboxes}
\end{figure}
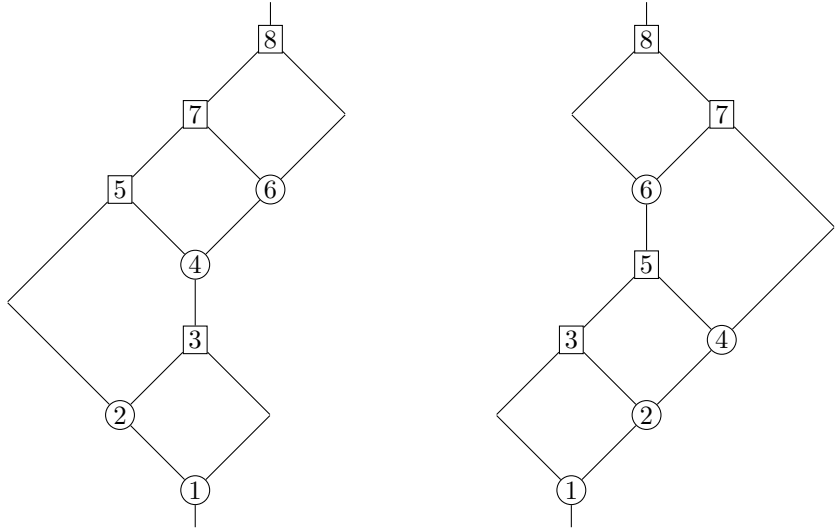

Figure~\ref{labelboxes} displays the labeling of operators for a
prograph having $4$ products and coproducts, and its reverse.

\begin{proposition}~\label{progtodyck}
  After having labeled product-coproduct prographs with the depth-left first
  algorithm, we can associate a step $(1,1)$ with values labeling
  coproducts and a step $(1,-1)$ with values labeling products. With such
  substitutions, the word $123\dots(2n)$ becomes a Dyck path.
  \begin{proof}
    The current height of the path is the number of active outputs
    minus one as the prograph is partially filled. Therefore the
    path remains over the horizontal axis. A primitive prograph (in
    the sense of an Hopf algebra element) is a prograph whose
    associated Dyck path returns to the horizontal axis only at the end.
  \end{proof}
\end{proposition}

\begin{definition}
  Let $n$ be a non-negative integer, we define a map $dw$ from
  prographs $PC(n)$ to weighted Dyck paths of length $2n$. We scan
  the prograph with depth-left first search labeling the operator from
  $1$ to $2n$ and starting a Dyck path at $(0,0)$ and reading the operator labeled
  by $i$ we build the weigthed Dyck path with the following rules.
  \begin{itemize}
  \item If $i$ labels a coproduct, we count the number $d$ of open
    outputs left free to the left of the grafting position of coproduct
    $i$. We add a step $(1,1)$ at the end of the Dyck path and we
    label this step with the integer $d$.
  \item If $i$ labels a product, we count the number $e$ of open
    outputs left free to the right of the grafting position of
    product $i$ (right from the right input of product labeled by
    $i$). We add a step $(1,-1)$ at the end of the Dyck path and we
    label this step with the integer~$e$.
  \end{itemize}
\end{definition}


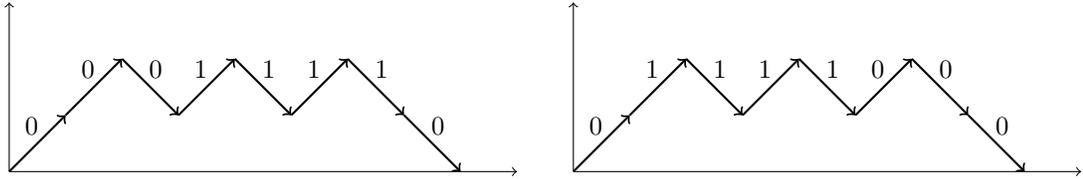
\begin{figure}
  \centering
  \begin{tikzpicture}[scale=0.75]
    \tikzstyle{p}=[->, thick]
    
    \draw[->] (0,0) -- (9, 0);
    \draw[->] (0,0) -- (0, 3);
    
    \draw[p] (0,0) -- (1,1);
    \draw[p] (1,1) -- (2,2);
    \draw[p] (2,2) -- (3,1);
    \draw[p] (3,1) -- (4,2);
    \draw[p] (4,2) -- (5,1);
    \draw[p] (5,1) -- (6,2);
    \draw[p] (6,2) -- (7,1);
    \draw[p] (7,1) -- (8,0);

    \draw (0.4, 0.5) node[above] {$0$};
    \draw (1.4, 1.5) node[above] {$0$};
    \draw (2.6, 1.5) node[above] {$0$};
    \draw (3.4, 1.5) node[above] {$1$};    
    \draw (4.6, 1.5) node[above] {$1$};
    \draw (5.4, 1.5) node[above] {$1$};
    \draw (6.6, 1.5) node[above] {$1$};
    \draw (7.6, 0.5) node[above] {$0$};

    \draw[->] (10,0) -- (19, 0);
    \draw[->] (10,0) -- (10, 3);
    
    \draw[p] (10,0) -- (11,1);
    \draw[p] (11,1) -- (12,2);
    \draw[p] (12,2) -- (13,1);
    \draw[p] (13,1) -- (14,2);
    \draw[p] (14,2) -- (15,1);
    \draw[p] (15,1) -- (16,2);
    \draw[p] (16,2) -- (17,1);
    \draw[p] (17,1) -- (18,0);

    \draw (10.4, 0.5) node[above] {$0$};
    \draw (11.4, 1.5) node[above] {$1$};
    \draw (12.6, 1.5) node[above] {$1$};
    \draw (13.4, 1.5) node[above] {$1$};    
    \draw (14.6, 1.5) node[above] {$1$};
    \draw (15.4, 1.5) node[above] {$0$};
    \draw (16.6, 1.5) node[above] {$0$};
    \draw (17.6, 0.5) node[above] {$0$};
    
  \end{tikzpicture}
  \caption{Weighted Dyck paths associated with the prographs of Figure~\ref{labelboxes}.}~\label{pathsweights}
\end{figure}

\begin{proposition}~\label{prop_weight}
  The map $dw$ is a bijection from $PC(n)$ to weighted Dyck
  paths of length $2n$ whose weight satisfies the following assertions:
  \begin{itemize}
    \item All weights of step $(1,1)$ are non-negative integers smaller than
      or equal to the starting height.
    \item All weigths of step $(1,-1)$ are non-negative integer smaller than
      or equal to the ending height.
    \item Weight are non decreasing on successive rises.
    \item Weight are non increasing on successive descents.
    \item On peaks of height $h$ where $d$ is the weight of a step
      $(1,1)$ just followed by a descent $(1,-1)$ labeled by $e$, we
      have: $e + d \leqslant h$.
    \item On valleys at height $h$ where $e$ is the weight of a step
      $(1,-1)$ just followed by a rise $(1,1)$ labeled by $d$, we
      have: $d + e \geqslant h$.
  \end{itemize}
\end{proposition}

Since our weighted Dyck paths in Proposition~\ref{prop_weight} are
Laguerre histories, we try in a first approach a customization of the
Françon-Viennot bijection~\cite{francon_viennot} to obtain up-down
permutations. All the variants we tested give up-down permutations which
do not necessarily avoid $(1234)$, therefore we need more to
solve Conjecture~\ref{conj1234}.

\begin{figure}
  \centering
  \scalebox{0.61}{
  \begin{tabular}{c|c|c|c|c}
  \begin{tikzpicture}[scale=0.55]
  \tikzstyle{prod}=[fill,draw,rectangle,minimum size=5pt,inner sep=1pt]
  \tikzstyle{cop}=[fill,draw,circle,minimum size=6pt,inner sep=1pt]

  \draw (-2,4) node[prod] (p3) {};
  \draw (-1,3) node[prod] (p2) {};
  \draw (0,2) node[prod] (p1) {};
  \draw (-2,2) node[cop] (c3) {};
  \draw (-1,1) node[cop] (c2) {};
  \draw (0,0) node[cop] (c1) {};

  \draw (0,0) -- (-3,3);
  \draw (0,0) -- (1,1);
  \draw (1,1) -- (-2,4);
  \draw (-3,3) -- (-2,4);
  \draw (-1,3) -- (-2,2);
  \draw (0,2) -- (-1,1);
  \end{tikzpicture}
  &
  \begin{tikzpicture}[scale=0.55]
  \tikzstyle{prod}=[fill,draw,rectangle,minimum size=5pt,inner sep=1pt]
  \tikzstyle{cop}=[fill,draw,circle,minimum size=6pt,inner sep=1pt]

  \draw (-1,5) node[prod] (p3) {};
  \draw (0,4) node[prod] (p2) {};
  \draw (-1,3) node[prod] (p1) {};
  \draw (-2,2) node[cop] (c3) {};
  \draw (-1,1) node[cop] (c2) {};
  \draw (0,0) node[cop] (c1) {};

  \draw (0,0) -- (-3,3);
  \draw (-2,2) -- (0,4);
  \draw (-1,1) -- (0,2);
  \draw (0,2) -- (-1,3);
  \draw (0,0) -- (2,2);
  \draw (2,2) -- (-1,5);
  \draw (-3,3) -- (-1,5);
  \end{tikzpicture}
  &
  \begin{tikzpicture}[scale=0.55]
  \tikzstyle{prod}=[fill,draw,rectangle,minimum size=5pt,inner sep=1pt]
  \tikzstyle{cop}=[fill,draw,circle,minimum size=6pt,inner sep=1pt]

  \draw (-1,5) node[prod] (p3) {};
  \draw (-2,4) node[prod] (p2) {};
  \draw (-1,3) node[prod] (p1) {};
  \draw (-2,2) node[cop] (c3) {};
  \draw (-1,1) node[cop] (c2) {};
  \draw (0,0) node[cop] (c1) {};

  \draw (0,0) -- (-3,3);
  \draw (-2,2) -- (-1,3);
  \draw (-1,1) -- (0,2);
  \draw (0,2) -- (-2,4);
  \draw (0,0) -- (2,2);
  \draw (2,2) -- (-1,5);
  \draw (-3,3) -- (-1,5);
  \end{tikzpicture}
  &
  \begin{tikzpicture}[scale=0.55]
  \tikzstyle{prod}=[fill,draw,rectangle,minimum size=5pt,inner sep=1pt]
  \tikzstyle{cop}=[fill,draw,circle,minimum size=6pt,inner sep=1pt]

  \draw (-1,5) node[prod] (p3) {};
  \draw (-2,4) node[prod] (p2) {};
  \draw (1,3) node[prod] (p1) {};
  \draw (-2,2) node[cop] (c3) {};
  \draw (-1,1) node[cop] (c2) {};
  \draw (0,0) node[cop] (c1) {};

  \draw (0,0) -- (-3,3);
  \draw (-2,2) -- (-1,3);
  \draw (-1,1) -- (1,3);
  \draw (-1,3) -- (-2,4);
  \draw (0,0) -- (2,2);
  \draw (2,2) -- (-1,5);
  \draw (-3,3) -- (-1,5);
  \end{tikzpicture}
  &
  \begin{tikzpicture}[scale=0.55]
  \tikzstyle{prod}=[fill,draw,rectangle,minimum size=5pt,inner sep=1pt]
  \tikzstyle{cop}=[fill,draw,circle,minimum size=6pt,inner sep=1pt]

  \draw (0,6) node[prod] (p3) {};
  \draw (-1,5) node[prod] (p2) {};
  \draw (-2,4) node[prod] (p1) {};
  \draw (-2,2) node[cop] (c3) {};
  \draw (-1,1) node[cop] (c2) {};
  \draw (0,0) node[cop] (c1) {};

  \draw (0,0) -- (-3,3);
  \draw (-3,3) -- (0,6);
  \draw (0,0) -- (3,3);
  \draw (3,3) -- (0,6);
  \draw (-1,1) -- (1,3);
  \draw (1,3) -- (-1,5);
  \draw (-2,2) -- (-1,3);
  \draw (-1,3) -- (-2,4);
  \end{tikzpicture}
  \\
  3 6 2 5 1 4
  &
  3 4 2 6 1 5
  &
  3 $\overline{5}$ 2 $\overline{6}$ 1 $\overline{4}$
  &
  3 6 2 4 1 5
  &
  3 5 2 4 1 6
  \\ \hline
  \begin{tikzpicture}[scale=0.55]
  \tikzstyle{prod}=[fill,draw,rectangle,minimum size=5pt,inner sep=1pt]
  \tikzstyle{cop}=[fill,draw,circle,minimum size=6pt,inner sep=1pt]

  \draw (-1,5) node[prod] (p3) {};
  \draw (0,4) node[prod] (p2) {};
  \draw (1,3) node[prod] (p1) {};
  \draw (0,2) node[cop] (c3) {};
  \draw (-1,1) node[cop] (c2) {};
  \draw (0,0) node[cop] (c1) {};

  \draw (0,0) -- (-3,3);
  \draw (-3,3) -- (-1,5);
  \draw (0,0) -- (2,2);
  \draw (2,2) -- (-1,5);
  \draw (-1,1) -- (1,3);
  \draw (0,2) -- (-1,3);
  \draw (-1,3) -- (0,4);
  \end{tikzpicture}
  &
  \begin{tikzpicture}[scale=0.55]
  \tikzstyle{prod}=[fill,draw,rectangle,minimum size=5pt,inner sep=1pt]
  \tikzstyle{cop}=[fill,draw,circle,minimum size=6pt,inner sep=1pt]

  \draw (0,6) node[prod] (p3) {};
  \draw (1,5) node[prod] (p2) {};
  \draw (0,4) node[prod] (p1) {};
  \draw (0,2) node[cop] (c3) {};
  \draw (-1,1) node[cop] (c2) {};
  \draw (0,0) node[cop] (c1) {};

  \draw (0,0) -- (-3,3);
  \draw (-3,3) -- (0,6);
  \draw (0,0) -- (3,3);
  \draw (3,3) -- (0,6);  
  \draw (-1,1) -- (1,3);
  \draw (0,2) -- (-1,3);
  \draw (-1,3) -- (1,5);
  \draw (1,3) -- (0,4);
  \end{tikzpicture}
  &
  \begin{tikzpicture}[scale=0.55]
  \tikzstyle{prod}=[fill,draw,rectangle,minimum size=5pt,inner sep=1pt]
  \tikzstyle{cop}=[fill,draw,circle,minimum size=6pt,inner sep=1pt]

  \draw (0,6) node[prod] (p3) {};
  \draw (-1,5) node[prod] (p2) {};
  \draw (0,4) node[prod] (p1) {};
  \draw (0,2) node[cop] (c3) {};
  \draw (-1,1) node[cop] (c2) {};
  \draw (0,0) node[cop] (c1) {};

  \draw (0,0) -- (-3,3);
  \draw (-3,3) -- (0,6);
  \draw (0,0) -- (3,3);
  \draw (3,3) -- (0,6);  
  \draw (-1,1) -- (1,3);
  \draw (0,2) -- (-1,3);
  \draw (-1,3) -- (0,4);
  \draw (1,3) -- (-1,5);
  \end{tikzpicture}
  &
  \begin{tikzpicture}[scale=0.55]
  \tikzstyle{prod}=[fill,draw,rectangle,minimum size=5pt,inner sep=1pt]
  \tikzstyle{cop}=[fill,draw,circle,minimum size=6pt,inner sep=1pt]

  \draw (0,4) node[prod] (p3) {};
  \draw (1,3) node[prod] (p2) {};
  \draw (-1,3) node[prod] (p1) {};
  \draw (0,2) node[cop] (c3) {};
  \draw (-1,1) node[cop] (c2) {};
  \draw (0,0) node[cop] (c1) {};

  \draw (0,0) -- (-2,2);
  \draw (-2,2) -- (0,4);
  \draw (0,0) -- (2,2);
  \draw (2,2) -- (0,4);
  \draw (-1,1) -- (1,3);
  \draw (0,2) -- (-1,3);
  \end{tikzpicture}
  &
  \begin{tikzpicture}[scale=0.55]
  \tikzstyle{prod}=[fill,draw,rectangle,minimum size=5pt,inner sep=1pt]
  \tikzstyle{cop}=[fill,draw,circle,minimum size=6pt,inner sep=1pt]

  \draw (1,5) node[prod] (p3) {};
  \draw (0,4) node[prod] (p2) {};
  \draw (-1,3) node[prod] (p1) {};
  \draw (0,2) node[cop] (c3) {};
  \draw (-1,1) node[cop] (c2) {};
  \draw (0,0) node[cop] (c1) {};

  \draw (0,0) -- (-2,2);
  \draw (-2,2) -- (1,5);
  \draw (0,0) -- (3,3);
  \draw (3,3) -- (1,5);
  \draw (-1,1) -- (1,3);
  \draw (0,2) -- (-1,3);
  \draw (1,3) -- (0,4);
  \end{tikzpicture}
  \\
  \underline{2} 6 \underline{1} 5 \underline{3} 4
  &
  \underline{2} 4 \underline{1} 6 \underline{3} 5
  &
  \underline{2} $\overline{5}$ \underline{1} $\overline{6}$ \underline{3} $\overline{4}$
  &
  \underline{2} 6 \underline{1} 4 \underline{3} 5
  &
  \underline{2} 5 \underline{1} 4 \underline{3} 6
  \\ \hline
  \begin{tikzpicture}[scale=0.55]
  \tikzstyle{prod}=[fill,draw,rectangle,minimum size=5pt,inner sep=1pt]
  \tikzstyle{cop}=[fill,draw,circle,minimum size=6pt,inner sep=1pt]

  \draw (-1,5) node[prod] (p3) {};
  \draw (0,4) node[prod] (p2) {};
  \draw (1,3) node[prod] (p1) {};
  \draw (0,2) node[cop] (c3) {};
  \draw (1,1) node[cop] (c2) {};
  \draw (0,0) node[cop] (c1) {};

  \draw (0,0) -- (-3,3);
  \draw (-3,3) -- (-1,5);
  \draw (0,0) -- (2,2);
  \draw (2,2) -- (-1,5);
  \draw (0,2) -- (1,3);
  \draw (1,1) -- (-1,3);
  \draw (-1,3) -- (0,4);
  \end{tikzpicture}
  &
  \begin{tikzpicture}[scale=0.55]
  \tikzstyle{prod}=[fill,draw,rectangle,minimum size=5pt,inner sep=1pt]
  \tikzstyle{cop}=[fill,draw,circle,minimum size=6pt,inner sep=1pt]

  \draw (0,6) node[prod] (p3) {};
  \draw (1,5) node[prod] (p2) {};
  \draw (0,4) node[prod] (p1) {};
  \draw (0,2) node[cop] (c3) {};
  \draw (1,1) node[cop] (c2) {};
  \draw (0,0) node[cop] (c1) {};

  \draw (0,0) -- (-3,3);
  \draw (-3,3) -- (0,6);
  \draw (0,0) -- (3,3);
  \draw (3,3) -- (0,6);
  \draw (1,1) -- (-1,3);
  \draw (-1,3) -- (1,5);
  \draw (0,2) -- (1,3);
  \draw (1,3) -- (0,4);
  \end{tikzpicture}
  &
  \begin{tikzpicture}[scale=0.55]
  \tikzstyle{prod}=[fill,draw,rectangle,minimum size=5pt,inner sep=1pt]
  \tikzstyle{cop}=[fill,draw,circle,minimum size=6pt,inner sep=1pt]

  \draw (0,6) node[prod] (p3) {};
  \draw (-1,5) node[prod] (p2) {};
  \draw (0,4) node[prod] (p1) {};
  \draw (0,2) node[cop] (c3) {};
  \draw (1,1) node[cop] (c2) {};
  \draw (0,0) node[cop] (c1) {};

  \draw (0,0) -- (-3,3);
  \draw (-3,3) -- (0,6);
  \draw (0,0) -- (3,3);
  \draw (3,3) -- (0,6);
  \draw (1,1) -- (-1,3);
  \draw (-1,3) -- (0,4);
  \draw (0,2) -- (1,3);
  \draw (1,3) -- (-1,5);
  \end{tikzpicture}
  &
  \begin{tikzpicture}[scale=0.55]
  \tikzstyle{prod}=[fill,draw,rectangle,minimum size=5pt,inner sep=1pt]
  \tikzstyle{cop}=[fill,draw,circle,minimum size=6pt,inner sep=1pt]

  \draw (0,4) node[prod] (p3) {};
  \draw (1,3) node[prod] (p2) {};
  \draw (-1,3) node[prod] (p1) {};
  \draw (0,2) node[cop] (c3) {};
  \draw (1,1) node[cop] (c2) {};
  \draw (0,0) node[cop] (c1) {};

  \draw (0,0) -- (-2,2);
  \draw (-2,2) -- (0,4);
  \draw (0,0) -- (2,2);
  \draw (2,2) -- (0,4);
  \draw (1,1) -- (-1,3);
  \draw (0,2) -- (1,3);
  \end{tikzpicture}
  &
  \begin{tikzpicture}[scale=0.55]
  \tikzstyle{prod}=[fill,draw,rectangle,minimum size=5pt,inner sep=1pt]
  \tikzstyle{cop}=[fill,draw,circle,minimum size=6pt,inner sep=1pt]

  \draw (1,5) node[prod] (p3) {};
  \draw (0,4) node[prod] (p2) {};
  \draw (-1,3) node[prod] (p1) {};
  \draw (0,2) node[cop] (c3) {};
  \draw (1,1) node[cop] (c2) {};
  \draw (0,0) node[cop] (c1) {};

  \draw (0,0) -- (-2,2);
  \draw (-2,2) -- (1,5);
  \draw (0,0) -- (3,3);
  \draw (3,3) -- (1,5);
  \draw (1,1) -- (-1,3);
  \draw (0,2) -- (1,3);
  \draw (1,3) -- (0,4);
  \end{tikzpicture}
  \\
  3 6 1 5 2 4
  &
  3 4 1 6 2 5
  &
  3 $\overline{5}$ 1 $\overline{6}$ 2 $\overline{4}$
  &
  3 6 1 4 2 5
  &
  3 5 1 4 2 6
  \\ \hline
  \begin{tikzpicture}[scale=0.55]
  \tikzstyle{prod}=[fill,draw,rectangle,minimum size=5pt,inner sep=1pt]
  \tikzstyle{cop}=[fill,draw,circle,minimum size=6pt,inner sep=1pt]

  \draw (-1,5) node[prod] (p3) {};
  \draw (0,4) node[prod] (p2) {};
  \draw (1,3) node[prod] (p1) {};
  \draw (-2,2) node[cop] (c3) {};
  \draw (1,1) node[cop] (c2) {};
  \draw (0,0) node[cop] (c1) {};

  \draw (0,0) -- (-3,3);
  \draw (-3,3) -- (-1,5);
  \draw (0,0) -- (2,2);
  \draw (2,2) -- (-1,5);
  \draw (0,2) -- (1,3);
  \draw (1,1) -- (0,2);
  \draw (-2,2) -- (0,4);
  \end{tikzpicture}
  &
  \begin{tikzpicture}[scale=0.55]
  \tikzstyle{prod}=[fill,draw,rectangle,minimum size=5pt,inner sep=1pt]
  \tikzstyle{cop}=[fill,draw,circle,minimum size=6pt,inner sep=1pt]

  \draw (0,4) node[prod] (p3) {};
  \draw (1,3) node[prod] (p2) {};
  \draw (0,2) node[prod] (p1) {};
  \draw (1,1) node[cop] (c3) {};
  \draw (-1,1) node[cop] (c2) {};
  \draw (0,0) node[cop] (c1) {};

  \draw (0,0) -- (-2,2);
  \draw (-2,2) -- (0,4);
  \draw (0,0) -- (2,2);
  \draw (2,2) -- (0,4);
  \draw (-1,1) -- (1,3);
  \draw (1,1) -- (0,2);
  \end{tikzpicture}
  &
  \begin{tikzpicture}[scale=0.55]
  \tikzstyle{prod}=[fill,draw,rectangle,minimum size=5pt,inner sep=1pt]
  \tikzstyle{cop}=[fill,draw,circle,minimum size=6pt,inner sep=1pt]

  \draw (0,4) node[prod] (p3) {};
  \draw (-1,3) node[prod] (p2) {};
  \draw (0,2) node[prod] (p1) {};
  \draw (1,1) node[cop] (c3) {};
  \draw (-1,1) node[cop] (c2) {};
  \draw (0,0) node[cop] (c1) {};

  \draw (0,0) -- (-2,2);
  \draw (-2,2) -- (0,4);
  \draw (0,0) -- (2,2);
  \draw (2,2) -- (0,4);
  \draw (-1,1) -- (0,2);
  \draw (1,1) -- (-1,3);
  \end{tikzpicture}
  \\
  2 6 3 5 1 4
  &
  2 4 3 6 1 5
  &
  2 $\overline{5}$ 3 $\overline{6}$ 1 $\overline{4}$
  \\ \cline{0-2}
  \begin{tikzpicture}[scale=0.55]
  \tikzstyle{prod}=[fill,draw,rectangle,minimum size=5pt,inner sep=1pt]
  \tikzstyle{cop}=[fill,draw,circle,minimum size=6pt,inner sep=1pt]

  \draw (0,6) node[prod] (p3) {};
  \draw (1,5) node[prod] (p2) {};
  \draw (2,4) node[prod] (p1) {};
  \draw (2,2) node[cop] (c3) {};
  \draw (1,1) node[cop] (c2) {};
  \draw (0,0) node[cop] (c1) {};

  \draw (0,0) -- (-3,3);
  \draw (-3,3) -- (0,6);
  \draw (0,0) -- (3,3);
  \draw (3,3) -- (0,6);
  \draw (1,1) -- (-1,3);
  \draw (-1,3) -- (1,5);
  \draw (2,2) -- (1,3);
  \draw (1,3) -- (2,4);
  \end{tikzpicture}
  &
  \begin{tikzpicture}[scale=0.55]
  \tikzstyle{prod}=[fill,draw,rectangle,minimum size=5pt,inner sep=1pt]
  \tikzstyle{cop}=[fill,draw,circle,minimum size=6pt,inner sep=1pt]

  \draw (1,5) node[prod] (p3) {};
  \draw (2,4) node[prod] (p2) {};
  \draw (1,3) node[prod] (p1) {};
  \draw (2,2) node[cop] (c3) {};
  \draw (1,1) node[cop] (c2) {};
  \draw (0,0) node[cop] (c1) {};

  \draw (0,0) -- (-2,2);
  \draw (-2,2) -- (1,5);
  \draw (0,0) -- (3,3);
  \draw (3,3) -- (1,5);
  \draw (1,1) -- (0,2);
  \draw (0,2) -- (2,4);
  \draw (2,2) -- (1,3);
  \end{tikzpicture}
  &
  \begin{tikzpicture}[scale=0.55]
  \tikzstyle{prod}=[fill,draw,rectangle,minimum size=5pt,inner sep=1pt]
  \tikzstyle{cop}=[fill,draw,circle,minimum size=6pt,inner sep=1pt]

  \draw (1,5) node[prod] (p3) {};
  \draw (0,4) node[prod] (p2) {};
  \draw (1,3) node[prod] (p1) {};
  \draw (2,2) node[cop] (c3) {};
  \draw (1,1) node[cop] (c2) {};
  \draw (0,0) node[cop] (c1) {};

  \draw (0,0) -- (-2,2);
  \draw (-2,2) -- (1,5);
  \draw (0,0) -- (3,3);
  \draw (3,3) -- (1,5);
  \draw (1,1) -- (0,2);
  \draw (0,2) -- (1,3);
  \draw (2,2) -- (0,4);
  \end{tikzpicture}
  \\
  1 6 3 5 2 4
  &
  1 4 3 6 2 5
  &
  1 $\overline{5}$ 3 $\overline{6}$ 2 $\overline{4}$
  \\
  \end{tabular}
  }
  \caption{The $21$ prographs associated with their up-down permutations
    of size $6$ avoiding $(1234)$ such that coproducts are labeled by $1,
    2$ and $3$ and products are labeled by $4, 5$ and $6$
    (\textit{eq.} $Vals(\sigma) = \{1,2,3\}$ and $Peaks(\sigma) =
    \{4,5,6\}$ on permutations).}~\label{matrix_pro}
\end{figure}
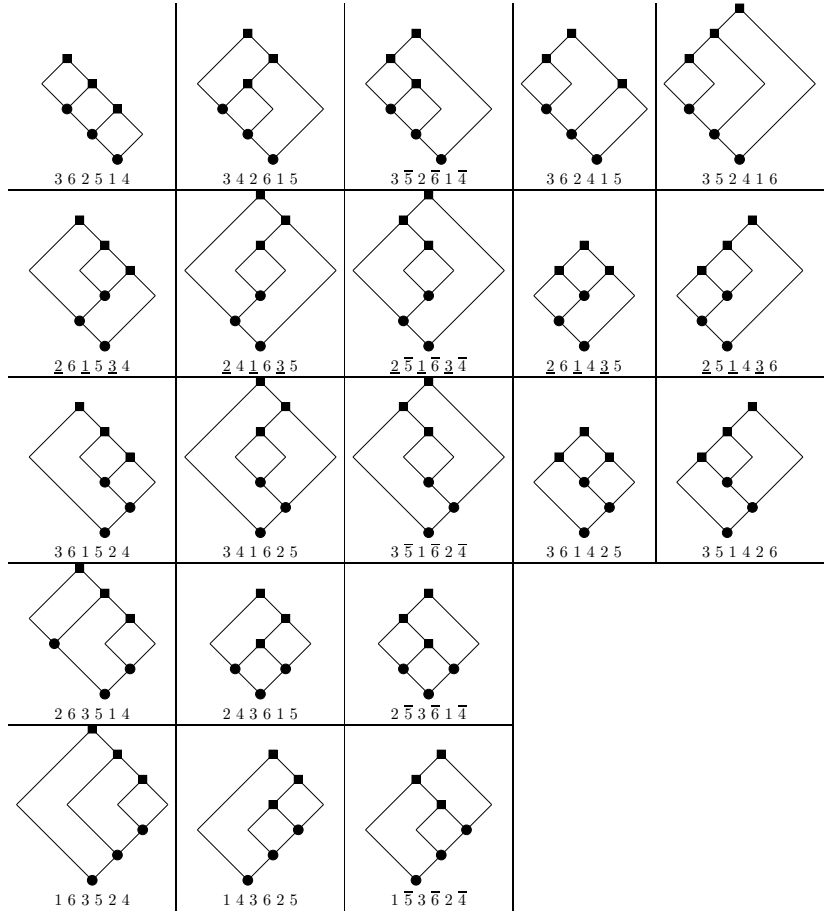

Nevertheless, we have a bijection for prographs of $PC(n)$ such that $Vals =
\{1,2, \dots n\}$ and $Peaks = \{n+1, n+2, \dots , 2n\}$ by using 
twice the bijection displayed in Figure~\ref{bij123_tree} and 
Proposition~\ref{prop_cond} applied with these special
conditions. Figure~\ref{matrix_pro} presents this bijection in size
$3$.

We are currently working on primitive prographs whose associated
weighted Dyck paths contain more than one peak. For all partitions of
$\{1,2, \dots ,2n\}$ into two sets $V$ and $P$, computer exploration
shows that the number of up-down permutations having for valleys $V$
and peaks $P$ is equal to the number of prographs whose labels of
coproducts are $V$ and labels of products are $P$. Therefore, we hope
to extend our bijection on product-coproduct prographs
making this new combinatorial class central for the study of objects
counted by the $3$-dimensional Catalan numbers.


The author thanks Samuele Giraudo for useful discussions and
comments. His experience with operads and PROs theories were
important to advise the author. This research was driven by computer
exploration using the open-source mathematical software
\texttt{Sage}~\cite{sage} and its algebraic combinatorics features
developed by the \texttt{Sage-Combinat} community~\cite{Sage-Combinat}.

\bibliographystyle{plain}
\bibliography{main}

\end{document}